\documentclass[11pt,a4paper]{article}

\title{Ring structure, uniform expressions and intersection homology}

\author{Jonathan Fine\relax
\thanks{203 Coldhams Lane, Cambridge, CB1 3HY, England.
\quad E-mail: \texttt{j.fine@pmms.cam.ac.uk}}
}
\date{11 May 1998}

\newcommand\PDelta{{\bf P\!}_\Delta}

\newcommand\bfR{{\bf R}}
\newcommand\bfZ{{\bf Z}}

\newcommand\vol{\mathop{\rm vol}\nolimits}

\newcommand\bibrule{\rule{2pc}{0.4pt}}

\begin{document}

\maketitle

\begin{abstract}
\noindent
Although intersection homology lacks a ring structure, certain expressions
(called uniform) in the intersection homology of an irreducible projective
variety $X$ always give the same value, when computed via the
decomposition theorem on any resolution $X_r\to X$.  This paper uses
uniform (and non-uniform) expressions to define what is believed to be the
usual intersection homology (and its local-global variant) of a convex
polytope (or a projective toric variety).  Such expressions are generated
by the facets, and so may lead to necessary numerical conditions on the
flag vector.  Most of the concepts, however, apply to more general
algebraic varieties, and perhaps some other situations also.
\end{abstract}

\section {Introduction}

This is the third of a series of papers devoted to intersection homology
and the combinatorics of convex polytopes.  The first \cite{bib.JF.LGIH}
gave a topological definition of the new local-global intersection homology
groups, and also gave the expected value of their Betti numbers, for
projective toric varieties (or equivalently convex polytopes).  The second
\cite{bib.JF.CPLA} constructed for each such group a complex which, if
exact at all but a single specified location, would produce a vector space
whose dimension is the expected value of the Betti number.  More briefly,
the paper gave both the Betti numbers and the topological definition,
while the second gave the linear algebra construction for the homology
groups.  As proofs of these statements have not been sought, they are at
present only conjectures.

The paper in some sense gives the geometric meaning of local-global
homology, at least in the context of convex polytopes.  It takes the view
that homology arises as the result of studying the calculation of
something, and of the properties of such a calculation.  Historically,
homology arose out of the study of intersection numbers (or, if one
prefers, the Schubert calculus), and was then placed by Poincar\'e into
its now well-known topological form.  This approach, however, does not
tell us what is the something that homology is calculating.  This paper
presents homology as a consequence of the study of volume, and its
dependence on the conditions defining the object being studied.  Such an
approach works well for convex polytopes but not, as noted in the final
section, for more general varieties.

This paper, together with its two companions, give the leading definitions
for the theory of local-global (intersection) homology, and also gives
conjectures for many of the main results that will hold, when this theory
is applied to convex polytopes.  (The major lack is any understanding of
how the various homology groups should be equipped with a ring-like
structure, that reduces to the usual one when the object is nonsingular.) 
The proofs are expected to be recursive, and involve a deeper knowledge of
the subject than we have at present.  The main purpose has been to present
what are believed to be the landmarks, and not to describe the paths
between them, nor to explore the intervening territory. Throughout convex
polytopes have been studied, because that are the least complicated
objects to which such a theory will apply.  By and large, this paper can
be read independently of the other two.

Convex polytopes satisfy subtle combinatorial inequalities, whose proof is
essentially geometric in nature.  The discovery of such results has been a
central motive for these papers.  (Other objects, such as graphs and
hypergraphs, will similarly satisfy subtle combinatorial inequalities, but
except for a few words in the final section, this matter is not taken up
again in this paper.  In \cite{bib.JF.QTHGFV,bib.JF.SFV,bib.JF.GFP} the
author introduces what he believes is the proper concept of a flag vector
for use in this new context.)

A central result in the combinatorics and geometry of convex polytopes is
the proof of the necessity (Stanley~\cite{bib.RS.NFSP}) and sufficiency
(Billera and Lee~\cite{bib.LB-CL.SMC}) of McMullen's conjectured
conditions~\cite{bib.McM.NFSP} on the face vector of simple polytopes. 
This paper studies general polytopes.  When $\Delta$ is a simple convex
polytope there is a homology ring $H_\bullet\Delta=\bigoplus_{i=0}^n
H_i\Delta$, whose properties imply McMullen's conditions.  The main goal of
this paper is to define a similar object for general polytopes.  This new
object will be based on the theory of intersection homology.  Although
formulated for convex polytopes, the definitions also apply to any
irreducible projective variety.

In some sense the central problem is this.  When $\Delta$ is simple the
homology object is a ring $H_\bullet\Delta$, generated over
$H_n\Delta\cong\bfR$ by $H_{n-1}\Delta$, or in other words the facets of
$\Delta$.  Although intersection homology does not supply a ring
structure, one wishes nonetheless to make a similar statement in the
general case.  The proposed solution to this problem is the concept of a
\emph{uniform expression}.  We will explain this first in the context of
algebraic varieties.

Suppose that $X$ is an irreducible projective variety of real dimension
$m=2n$.  Suppose also that $i+j=m$.  In that case (middle perversity)
intersection homology provides a pair of homology groups, which we will
denote by $H_i\Delta$ and $H_j\Delta$, together with the pairing
\[
    \frown \, : H_iX \otimes H_jX \to H_0X \cong \bfR
\]
that is, by Poincar\'e duality, nondegenerate.  If $X$ were nonsingular,
then for any $i$ and $j$ there would be a product map
\[
    \frown \, : H_iX \otimes H_jX \to H_kX
\]
that gives the homology ring structure.  Here $(m-k) = (m-i) + (m-j)$.

Now suppose that $X_1 \to X$ is a resolution of the singularities of $X$. 
According to the decomposition theorem \cite{bib.AB-JB-PD.FP}, there is
now an inclusion $H_iX \hookrightarrow H_iX_1$ of the (intersection)
homology groups of $X$ into the ordinary homology of $X_1$ (which is the
same as its intersection homology, because $X_1$ is nonsingular).  The
decomposition theorem is a deep result, and the exact nature of the
inclusion it provides is not well understood. Nonetheless, it exists, and
we can use it.  (In fact, we must assume that $X_1 \to X$ is a resolution
of $X$ together with what is known as a relatively ample line bundle.  But
this detail does not affect the exposition.)

We can now use the decomposition theorem, together with Poincar\'e
duality, to provide a ring-like structure on $H_\bullet X$.  First suppose
$i+j+k=m$, and that $\psi$, $\eta$ and $\xi$ are homology classes on $X$
with dimensions $i$, $j$ and $k$ respectively.  Now let
\[
    X_1 \frown \psi \frown \eta \frown \xi 
        \in H_0X_1 \cong H_0X \cong \bfR
\]
denote the result of lifting the classes to $X_1$, and then computing
their intersection number.  The value will of course depend on the choice
of $X_1 \to X$.

Now think of $\psi$ and $\eta$ as fixed, and $\psi$ as variable.  The
above expression is a linear function of $\psi$ and so, by Poincar\'e
duality, it can be represented by an intersection homology class in
$H_{m-k}X$, which we can denote by $\psi \frown_1 \eta $.  As noted, it
will in general depend on the choice of $X_1\to X$.  It does not seem
likely that such a product will be associative, although it may have some
quasi-associative properties.

Now say that an expression $\zeta = \psi \frown \eta$ (or more generally
$\zeta = \psi_1 \frown \dots \frown \psi_s$) in the (intersection)
homology classes of $X$ is \emph{uniform} if the evaluation of
\[
    X_r \frown \zeta \frown \xi \in H_0X_r \cong H_0X
\]
does not depend on the choice of a resolution $X_r \to X$, where $\xi$ is
an arbitrary cycle in $H_\bullet X$, of complementary dimension.  The same
concept will also be applied to sums of products of cycles, provided each
term in the expression has the same degree.  Clearly, it follows from
Poincar\'e duality that each uniform expression determines an
(intersection) homology class, provided we assume as we do that $X$ has at
least one resolution $X_r \to X$.

The concept of uniform expression can now be used to make precise the
statement that the homology object $H_\bullet\Delta$ of a general convex
polytope $\Delta$ is generated, as in the simple case, by $H_{n-1}\Delta$. 
In fact, because the inclusion provided by the decomposition theorem is in
general still rather obscure, we will reverse the process, and define the
(intersection) homology groups $H_i\Delta$ of a general $\Delta$ to be
certain uniform expressions, modulo the relations induced by evaluation
(or Poincar\'e duality).  This will be done in the next section.

Expressions that are not uniform can also be useful.  They are capable of
recording information about the singularities of $X$.  (Uniform
expressions are by definition incapable of doing this.)  Consider for
example the expression
\[
    (X_1 - X_2 ) \frown \psi \frown \eta \frown \xi \in H_0X \> ,
\]
where $X_1 \to X$ and $X_2 \to X$ are two distinct resolutions of $X$.  If
$\psi \frown \eta$ is not uniform then this expression can be non-zero. 
Its value will be concentrated, so to speak, along the region of $X$ at
which the two resolutions differ.  In \S3 a scheme will be presented that
records such information in what are believed to be variant forms of the
local-global intersection homology groups that the author has introduced
elsewhere~\cite{bib.JF.LGIH,bib.JF.CPLA}.

\section{Volume and homology}

There is a connection between the homology and the volume of a simple
convex polytope, that deserves to be more widely known.  It leads to a
presentation of the homology theory, equivalent to the usual
one~\cite{bib.VD.GTV}, which is close in spirit to the concept of uniform
expressions.

Throughout $\Delta$ will be an $n$-dimensional convex polytope in an
$n$-dimensional affine space, described by an irredundant system
\[
    \alpha_i v \geq 0 \qquad i = 1 , \ldots , f_{n-1}
\]
of affine linear inequalities, one for each facet $\delta_i$ of $\Delta$. 
This writes $\Delta$ as an intersection of half-spaces.  We now subject
the half-spaces to a small displacement, and consider how the volume
changes.  In other words, let $\Delta(\epsilon)$ be the convex polytope
\[
    \alpha_i v \geq \epsilon_i \qquad i = 1 , \ldots , f_{n-1}
\]
where the $\epsilon_i$ are close to zero.

This process of $\epsilon$-variation will in general change the
combinatorial structure of $\Delta$.  The cone (or pyramid) on a square is
an example of this.  It will be studied further, in the course of this
paper.  If the $\epsilon$-variation process does not change the
combinatorics (for $\epsilon$ close to zero of course) then $\Delta$ will
be said to be \emph{simple}.  This happens exactly when at each vertex
there are exactly $n$ facets (or equivalently $n$ edges), which is the
usual definition of simple.

Now suppose that $\Delta$ is a simple polytope.  From this it follows
without any real difficulty that the volume $\vol(\epsilon)$ of
$\Delta(\epsilon)$ is a polynomial of degree $n$ in the linear quantity
$\epsilon$, at least for $\epsilon$ small.  (To prove this, decompose
$\Delta(\epsilon)$ into pyramids, one for each facet.  The height of each
such pyramid will be linear in $\epsilon$, while by induction the area of
the base will be a polynomial of degree $n-1$ in $\epsilon$.  For general
polytopes the volume $\vol(\Delta)$ will be given by one of a number of
polynomials, one for each simple structure reached via small
$\epsilon$-variation.)

The polynomial $\vol(\Delta)$ is perhaps complicated, but it is far from
arbitrary.  It has a great deal of structure.  For example, to know the
top degree term $\vol_\Delta$ is equivalent to knowing the homology ring
$H_\bullet\Delta$, a fact that we will explain shortly.  The lower degree
terms give information about the Chern classes of $\Delta$ (or of the
projective toric variety $\PDelta$ associated to $\Delta$).  This we do not
need, and so will not discuss further.

The top degree term $\vol_\Delta$ can be written as a symmetric
multilinear form, also to be denoted by $\vol_\Delta$, in $n$ variables
$\alpha_1$, $\ldots$, $\alpha_n$.  Each $\alpha_i$ is a (small)
$\epsilon$-variation $\Delta$.  The simplest case is where each $\alpha_i$
is a displacement of a single facet of $\Delta$.  By linearity, these
determine the general case, which is where the $\alpha_i$ are formal sums
of displacements of individual facets.  In any case, let $T^1\Delta$
denote the linear span of all such $\alpha_i$.  It has dimension
$f_{n-1}$, the number of facets, and represents \emph{thickenings} or
$\epsilon$-displacements of the facets.

Now let $T^i\Delta$ denote the $i$-fold tensor product of $T^1\Delta$ with
itself.  As already noted, $\vol_\Delta$ can be thought of as a symmetric
linear function on $T^n\Delta$.  Now let $i+j=n$ and consider the map
\begin{equation}
\label{eqn.Ti-Tj}
    \vol_\Delta : T^i \Delta \otimes T^j \Delta \to \bfR
\end{equation}
induced by $T^i\Delta \otimes T^j \Delta \cong T^n \Delta$.  This is a
pairing of vector spaces.  Now form the null spaces
\[
    N^i\Delta = \{ \, \psi \in T^i 
                    \> | \> 
                    \vol_\Delta ( \psi , \eta ) = 0
                    \mbox{ for all $\psi \in T^j$} \, \}
\]
and $N^j\Delta$ similarly, and then form the quotient spaces $H^i =
T^i/N^i$, $H^j = T^j/N^j$.  The pairing now descends to give a
nondegenerate (or perfect) pairing
\begin{equation}
\label{eqn.Hi-Hj}
    H^i \Delta \otimes H^j \Delta \to \bfR
\end{equation}
between a pair of vector spaces.

This definition of the (co)homology spaces $H^i\Delta$ is equivalent to
the usual one.  (To simplify matters, we will identify the cohomology
group $H^i\Delta$ with the homology group $H_j\Delta$, where $i+j=n$.
Thus, cohomology is homology, but indexed by codimension as a superscript,
rather than dimension as a subscript.)  The key facts are this.  First,
the (co)homology is generated by the facets.  This means that $T^i\Delta$
will generate $H^i\Delta$.  Second, Poincar\'e duality holds.  This means
that passing from (\ref{eqn.Ti-Tj}) to (\ref{eqn.Hi-Hj}) will produce
(co)homology in the usual sense.  Third, the volume of $\Delta$ and the
volume of the toric variety $\PDelta$ agree, provided $\PDelta$ exists. 
More exactly, $\Delta$ determines on $\PDelta$ a `hyperplane class' or
K\"ahler form $\omega$ such that
\[
    \vol \PDelta = \int _{\PDelta} \omega ^n
\]
is proportional to $\vol\Delta$.  This is true for all $\Delta$, and so
respects $\epsilon$-variation.  Agreement here is enough to force
agreement everywhere.

One can look at this in another way.  Suppose $H^i\Delta$ is known, and
that it is generated by $H^1\Delta$.  The induced map
\begin{equation}
\label{eqn.H1-n-times}
    H^1\Delta \otimes \dots \otimes H^1 \Delta \to H^n\Delta \cong \bfR
\end{equation}
on the $n$-fold tensor product is more-or-less equivalent to $\vol_\Delta$
on $T^n\Delta$.  (In fact, $H^1\Delta$ is $T^1\Delta$ modulo certain
relations.)  Because of Poincar\'e duality (and generation by
$H^1\Delta$), using $i+j=n$ to break this up into a pairing will produce
$H^i$ and $H^j$.  The induced map \ref{eqn.H1-n-times} is however
equivalent to knowing a homogeneous polynomial on $H^1\Delta$ (or
$T^1\Delta$).  In this way, the $\epsilon$-variation volume polynomial and
the homology theory determine each other.

Here is a result that will be used to help do an example in the next
section.  Clearly, an $\epsilon$-variation that is in fact a rigid bodily
translation of $\Delta$ will not change its volume.  Thus,
$\vol_\Delta(\alpha_1, \dots , \alpha_n)$ will be zero if any one of the
$\alpha_i$ is of this form.  (Each $\alpha_i$ is a formal sum of facet
displacements.)  Next, just as
\begin{equation}
\label{eqn.2xy}
    2 \langle x | y \rangle
    =  \| x+y \| ^2 - \| x\| ^2 - \| y\|^2
\end{equation}
turns a quadratic form into a bilinear form, so a similar process will
produce the multilinear form $\vol_\Delta$ from the polynomial form.  It
now follows that if each $\alpha_i$ is a displacement of a single facet
$\delta_i$, and these facets $\delta_i$ have empty common intersection,
then $\vol(\alpha_i,\dots,\alpha_n)$ is zero.  For $n=2$, this is because
the $\delta_1$ and $\delta_2$ displacements do not, so to speak, interfere
with each other, and so (\ref{eqn.2xy}) will be zero.  A similar argument
holds in general.

In the usual approach Poincar\'e duality asserts that a pairing between
already defined homology spaces is nondegenerate.  In the present
approach, duality is part of the definition.  Poincar\'e duality then
becomes a characterization of the nullspaces $N^i$ and $N^j$.  It asserts
that they are generated by the rigid displacement and empty intersection
results stated in the previous paragraph.

Finally, in the present approach, the ring structure is automatic. 
Clearly there is a map $T^i\otimes T^j \to T^{i+j}$, and it is obvious that
$N^i\otimes T^j$ and $T^j\otimes N^j$ both lie in $N^{i+j}$.  Hence the
tensor product descends to give a multiplication $H^i\otimes H^j\to
H^{i+j}$ on (co)homology.  Because all this is based on the well-defined
concept of volume, there is no need to move cycles into general position
and so forth.

\section {Uniform expressions}

When $\Delta$ is not simple, the polytopes produced by
$\epsilon$-variation will in general have different combinatorial
structures.  The cone (or pyramid) on a square is the simplest example. 
Suppose that $\Delta$ is this polytope.  Let the points of the compass
$N$, $S$, $E$ and $W$ denote the triangular facets, and let $B$ denote the
square base.  Use the same symbols to denote thickenings of these facets,
so chosen that $N-S$, $E-W$ and $N+E-B$ each represent a thickening due to
a rigid displacement.  (Such displacements do not change the volume.)
Normalise so that $N\frown E \frown B$ has volume one.  Now let
$\Delta_{NS}$ denote the simple polytope that results from `squeezing in'
the $E$ and $W$ faces (or `letting out' the $N$ and $S$ faces).  Thus,
$\Delta_{NS}$ is a `ridge roof' polytope, with the ridge running $NS$. 
Similarly, let $\Delta_{EW}$ be the other ridge roof polytope obtained
from $\Delta$ by the $\epsilon$-variation process.

As before, let $N$, $S$, $E$, $W$ and $B$ denote the previous thickened
facets on $\Delta_{NS}$ and $\Delta_{EW}$, as well as on $\Delta$.  That
$N-S$, $E-W$ and $N+E-B$ represent trivial thickenings is still true.  In
the language of algebraic varieties, $\Delta_{NS}$ and$\Delta_{EW}$ are
two different resolutions of $\Delta$.  The $H_2$ groups of the three
polytopes are naturally isomorphic, and this isomorphism is the inclusion
given by the decomposition theorem.  Analogous results will hold for all
convex polytopes.  That $\epsilon$-variation does not change `the geometry
in codimension one' allows us to use facets as a starting point for the
use of the decomposition theorem.

On $\Delta$, there are five generators for $H_2\Delta$ (the facets) and
three relations ($N-S=E-W=N+E-B=0$).  This leaves a two dimensional space,
with generators $N$ and $E$ say.  We will now look for uniform
expressions.  In other words, adapting the notation of the previous
section, we seek $\alpha$ in $T^2\Delta$ such that
\[
    \Delta_{NS} \frown N \frown \alpha 
        = \Delta_{EW} \frown N \frown \alpha
\]
and so on for the other facets.  More exactly we seek an $\alpha$ and a
$\beta$, both with this uniform property, that provide a basis for $H_1$
that is dual to the $\{\,N,E\,\}$ basis for $H_2$.

Here is an easy solution.  Take $\alpha=W\frown B$ and $\beta=S\frown B$. 
In every case the intersection or homology calculation can be performed
entirely on the base $B$ of $\Delta$, and so is insensitive to the choice
of $\Delta_{NS}$ or $\Delta_{EW}$.  Thus, $\alpha$ and $\beta$ are uniform
expressions, that provide the sought for dual basis.

This solution is perhaps too easy.  Duality is in essence a local matter. 
About the apex of $\Delta$ both $N$ and $E$ represent non-trivial local
cycles.  However, because $N+E-B$ is trivial, and because $B$ is not local
to the apex, the local cycles due to $N$ and $E$ are equal but opposite. 
We therefore would like a uniform expression that is both local to the
apex and dual to $N$.

The quantity $\eta=E\frown W - N \frown S$ is a solution to this harder
local problem.  The verification proceeds as follows.  First
\[
    \Delta_{NS} \frown N \frown ( E\frown W - N \frown S )
    = \Delta_{NS} \frown N \frown E \frown W - 0 = 1
\]
because on $\Delta_{NS}$ the facets $N$ and $S$ do not meet.  Next
\[
    \Delta_{EW} \frown N \frown ( E\frown W - N \frown S )
    = 0 - \Delta_{EW} \frown N \frown N \frown S
\]
because on $\Delta_{EW}$ the facets $E$ and $W$ do not meet.  To compute
$N \frown N \frown S$ on $\Delta_{EW}$ use $N=B-E$ to obtain
\[
      \Delta_{EW} \frown B \frown N \frown S
    - \Delta_{EW} \frown E \frown N \frown S
    = 0 -1 = -1
\]
(because $B$, $N$ and $S$ have empty intersection on $\Delta_{EW}$) and
thus the equation
\[
    \Delta_{NS} \frown N \frown \eta 
    = \Delta_{EW} \frown N \frown \eta
\]
holds.  The verification that
\[
    \Delta_{NS} \frown E \frown \eta 
    = \Delta_{EW} \frown E \frown \eta
\]
is left to the reader.  (It also follows because $N+E=B$, and any
evaluation involving $B$ is automatically uniform.)

Now let $\Delta$ be an arbitrary convex polytope, as usual of dimension
$n$.  In the simple case there was a single (multilinear) volume form
$\vol_\Delta$, whose evaluation on $T^i\otimes T^j$ (with $i+j=n$) induced
the homology groups.  In the general case we will define $\vol_\Delta$ to
be the set $\{\,\vol_r\,\}$ of multilinear forms, due to \emph{all} the
simple polytopes $\Delta_r$ that arise from $\Delta$ as a result of
$\epsilon$-variation.  We now seek subspaces $U^i$ and $U^j$ of $T^i$ and
$T^j$, such that $\vol_\Delta$ is on $U^i\otimes U^j$ a single valued
function.

The task now is to define $U^i$ and $U^j$.  Of course, once say $U^i$ is
known, one can define $U^j$ to be all $\xi$ in $T^j$ such that
$\vol_\Delta(\eta,\xi)$ is single valued, for any $\eta$ in $U^i$.  One
can then as in the simple case quotient by the null spaces $N^i$ and $N^j$
of $U^i\otimes U^j \to \bfR$ to obtain a perfect pairing on $H^i$ and
$H^j$.

The correct choice of the $U^i$ and $U^j$ is somewhat delicate.  One would
like in some sense to choose them so that $H^i$ and $H^j$ are as large as
possible.  In this way, they will contain as much information as they are
able to.  Enlarging $U^i$, if it does not change $U^j$, might enlarge
$H^i$.  But it might reduce $U^j$, and thus enlarge $N^i \subseteq U^i$. 
As noted, once say $U^i$ is known, the rest of the construction follows in
a mechanical manner.  Thus, the definition of $U^i$ is the central
question of this section.

For $i=1$ the correct value of $U^1$ is already known.  It is the whole of
$T^1$.  Thus, $U^{n-1}$ consists of all $\xi$ in $T^{n-1}$ such that
$\vol_\Delta(\eta,\xi)$ is single valued, for any $\eta$ in $T^1$.  For
this construction, and indeed the whole of the uniform expression
programme, to be correct, the resulting homology groups should have the
desired dimension, which in this case is $f_{n-1}-n$.  In other words,
there should be enough uniform $\xi$.

In the present case of $U^{n-1}$ we can use the strong Lefschetz theorem
to find such $\xi$.  Choice of a point $p$ interior to $\Delta$ will
determine a thickening $\omega_p \in T^1\Delta$, namely the displacement
that would take each of the facets from $p$ to where they actually are. 
Changing $p$ will add a rigid (and so trivial) displacement to $\omega_p$. 
We will think of $\omega_p$ as the unique \emph{hyperplane class} or
\emph{Lefschetz element} $\omega$ in $T^1\Delta$, although it is only when
in $H^1\Delta$ that it becomes unique.

The Lefschetz element $\omega$ has an important local property, and an
important global property.  The local property is that it is, in the
language of algebraic geometry, a \emph{slice} or a \emph{section} through
a variety. As a result, it is locally trivial, in the sense that at each
vertex there is a (unique) trivial displacement that agrees with $\omega$
on the facets through that vertex.  We will let $S^1$ denote all such
locally trivial systems of thickened facets.  Similarly, $S^i$ will be the
$i$-fold tensor product of $T^i$.  If $\Delta$ is simple, then $S^1=T^1$,
and vice versa.  Even though in general two intersection homology classes
cannot be multiplied together, they always can be, provided one of them is
locally trivial.

The global property is the strong Lefschetz theorem.  This tells us that
for $i+j=n$ and $i \geq j$ the map
\[
    \omega ^ {j-i} : H^j \Delta \to H^i \Delta
\]
is an isomorphism.  This is at present known only for simple polytopes,
general polytopes with rational vertices, and some other special cases. 
But it is reasonable to assume that it is true in general.

As noted, $U^1$ is $T^1$.  Now let $U^{n-1}$ be the expressions (in
$T^{n-1}$) that are uniform when paired against $U^1$.  The first property
of $\omega$ allows us to conclude that $\omega^{n-2}\xi$ is a uniform
expression, for any $\xi \in U^1 = T^1$.  (This follows because $\Delta$
is simple along its faces of dimension $n-2$.)  The second property allows
us to conclude that there are enough such elements of $U^{n-1}$, to obtain
dual spaces $H^1$ and $H^{n-1}$ of dimension $f_{n-1}-n$.  (This requires
either Poincar\'e duality for intersection homology, or better some
further facts about the strong Lefschetz theorem.  The required facts are,
for $U^1$ and $U^{n-1}$, a consequence of Minkowski's results
\cite{bib.HM.ALKP} on facet area and `outward normal vectors'.)  Thus, it
is already known that at least in this case the uniform expression
approach gives the correct answer.  There is unlikely to be an easy proof
of this result.

We now turn to $U^2$ and $U^{n-2}$, or more generally $U^i$ and $U^j$, for
$i+j=n$ and $i \leq j$.  We can suppose that $U^{i-1}$ is already known. 
If $\xi$ is to lie in $U^i$, the product
\begin{equation}
\label{eqn.xi-omega-eta}
    \xi \frown \omega ^ k \frown \eta \in T^n
\end{equation}
must also be uniform, for any $\eta$ in $U^{i-1}$.  Here, $k$ is $n-2i+1$.
We will assume that this necessary condition is also sufficient, or in
other words that it produces the correct definition of $U^i$.

The space $U^j$ of complementary dimension can be defined just as before,
to be the expressions that are uniform when paired with expressions in
$U^i$.  This completes the definition except for the case $i=j$ (and $n$
is even).  In this case $U^i$ and $U^j$ has better be equal.  But this
question arises already, with the other $U^i$.  If $U^i$ has been
correctly defined then expressions such as (\ref{eqn.xi-omega-eta}) must
again be uniform, where now $\xi$ and $\eta$ are in $U^i$, and $k$ is
$n-2i$.  This then is a wished for property of $U^i$.  If it fails to
hold, then the uniform expression programme also fails.  The equality of
$U^i$ and $U^j$ for $i=j$ is a special case of this.  Finally, the
definitions of $U^0$ and $U^n$, and of $U^1$ and $U^{n-1}$, are special
cases of this general scheme.

\section {Local-global intersection homology}

By design, the unifrom expressions of the previous section are insensitive
to the choice of a resolution (simple $\epsilon$-variation) of the object
being studied.  The definition of this section will record information
about how the value of a (non-uniform) expression depends on the
resolution chosen.

Here is an example.  Let $\Delta$ be the cone on a square, and let
$\Delta_{NS}$, $\Delta_{EW}$, $N$, $S$, $E$, $W$ and $B$ be as before. 
Now consider the expression $N\frown S \frown E$.  This expression is not
uniform.  On $\Delta_{NS}$ it evaluates to zero, because the $N$ and $S$
facets do not meet.  On $\Delta_{EW}$ it is $1$.  (Replace $N$ by the
equivalent cycle $B-W$.  On $\Delta_{EW}$ the facets $W$ and $E$ do not
meet.  Moreover, $B\frown S \frown E$ is equal to $1$.)

The task is to organise such information in a sensible way.  The author
has already defined a theory of local-global intersection homology
\cite{bib.JF.LGIH,bib.JF.CPLA} which will, usually behind the scenes,
guide the definitions that follow.  This theory provides, for each
$n$-polytope, a system of $F_{n+2}$ `Betti numbers', of which $F_{n+1}$
are linearly independent.  These `Betti numbers' are organised into
$F_{n}$ sequences or strings.  The Lefschetz operator $\omega$ goes from
the homology group at one location in a string to the previous one (or
next, depending on the indexing scheme).  Here, $F_i$ is the $i$-th
Fibonacci number.  Bayer and Billera \cite{bib.MB-LB.gDS} showed that the
flag vectors of $n$-polytopes span a space whose dimension is $F_{n+1}$. 
The `Betti numbers' are a re-encoding of information in the flag vector.

Now let $\Delta$ be any $3$-polytope, and let $\Delta_1$ be any resolution
(by which we mean of course, a simple $\epsilon$-variation).  Now think of
the expression
\[
    \Delta_1 \frown \alpha \frown \beta \frown \gamma \qquad 
        \alpha, \beta, \gamma \in U^1 = T^1
\]
as a linear function of $\gamma$.  It will of course vanish for $\gamma
\in N^1$, and so by Poincar\'e duality determines a class in $H_1\Delta =
H^2\Delta$.  (In the first section, this class was denoted by $\alpha
\frown_1 \beta$.)

If $\Delta_1$ were replaced by another resolution $\Delta_2$, a different
class might result.  Now consider the expression
\[
    (\Delta_1 - \Delta_2 ) \frown \alpha \frown \beta \frown \gamma 
        \qquad \alpha, \beta, \gamma \in U^1 = T^1
\]
as a linear function of $\gamma$.  Again it represents a class in
$H_1\Delta$, but in contrast to the previous case it has an important
geometric property.  It is concentrated about, or local to, the locus in
$\Delta$ above which the two resolutions $\Delta_1$ and $\Delta_2$ differ. 
In the present case ($3$-polytopes) this locus will be a finite set of
points.  Thus, global cycles that are constructed in this way satisfy
certain locality properties.  This is why they are called local-global
cycles.  They are neither purely local or purely global, but partake in
the properties of both.

Now consider the octahedron.  Such a local-global cycle can be found at
each of its six vertices.  The octahedron has $(1,5,5,1)$ as its Betti
numbers, and so there must be, when considered as global cycles, a
relation between these six local-global cycles.  (In fact there are two
such relations.)  These relations do not however respect the local aspect
of the given local-global cycles.  It is important, particularly in higher
dimensions, to have a satisfactory definition of the equivalence of
local-global cycles.

The solution to this problem is to introduce a dual theory.  In other
words, a bilinear pairing will be produced.  As in the previous section
the quotients by the null-spaces will be called the \emph{local-global
homology groups}.  In contrast to the previous case, where $H^i$ and $H^j$
were paired with each other, here the pairing will be between
\emph{compact} and \emph{open} local-global cycles.  They are markedly
different objects. The cycles just recently introduced are compact.

We continue to study the case $n=3$.  As already noted, any global cycle
(thickened facet) $\gamma \in U^1 = T^1$ can be paired with a compact
local-global cycle such as
\begin{equation}
\label{eqn.Delta12-alpha-beta}
    (\Delta_1 -\Delta_2) \frown \alpha \frown \beta
\end{equation}
but that there will not in general be enough such cycles.  We can obtain
more by relaxing the properties that $\gamma$ satisfies.  This involves an
understanding of how intersection numbers (products of homology cycles)
can be calculated.

Suppose $\Delta$ is simple.  Previously the intersection number of $n$
thickened facets was defined using the already well-defined polynomial
$\vol\Delta(\epsilon)$.  There is another method, which has been used in
the examples.  It is to use equivalence of cycles to move them into
`general position', and to then compute the well-defined volume of their
common intersection (as thickened facets).  This is the traditional
method.  It does not rely on the volume polynomial.  That it produces a
well-defined product is a consequence of two facts.  The first is that
cycles can always be moved into general position.  This ensures that the
product can be calculated in every case.  The second fact is that
howsoever the calculation is performed, the same answer results.

The dual theory of \emph{open local-global cycles} will be developed using
this `general position' approach to the pairing.  When a product of $n$
thickened facets is in general position, the intersection takes place at a
vertex.  Suppose now that instead of having a global thickened facet
$\gamma$, one has at each vertex $v$ a thickened facet $\gamma_v$.  This
is a more general concept, for $\gamma$ at $v$ need not be the same
thickened facet as $\gamma$ at $w$.  Such will be called a (vertex
centered) system of open local cycles, as will be formal sums of such. 
Each (formal sum of) global thickened facets will determine such a system,
although in general the converse is plainly false.

Recall that an intersection pairing will exist as a consequence of the two
properties of computability and consistency.  The first is essentially a
local matter.  To move a class from a vertex, add to it a suitable cycle
that is equivalent to zero.  This process can be used to produce a perhaps
ill-defined pairing between compact cycles on $\Delta$ and vertex-centered
open local cycles of complementary dimension, as defined in the previous
paragraph.  The cycles equivalent to zero are as in the usual theory, but
understood as open local cycles.

In general, such a product will not be consistent.  However, we wish to
apply it not to all compact cycles, but only to those that arise in a
particular way, namely as compact local-global cycles.  Suppose such a set
of cycles is given.  For an open local cycle to give a consistent answer
when evaluated against such a family of compact cycles is a \emph{global}
property of open local cycles.  Such a cycle will be called, as one would
expect, an open local-global cycle.  The global cycles (of appropriate
dimension) have this property.  Beause the open cycles are being paired
with a restricted set of compact cycles, there will in general be many
open cycles that do not arise from a global cycle.  In fact, for $n=3$ the
consistency condition turns out to be vacuous.  (In higher dimensions it
is more subtle.)

There are some points that should be clarified.  First, when computing the
compact-open local-global pairing, any moving that is to be applied on the
compact side should respect the local-global origin of the cycles.  In
other words, a cycle such as (\ref{eqn.Delta12-alpha-beta}) should be
moved only by moving $\alpha$ and $\beta$, and not by the arbitrary moving
of elements in $H_1\Delta_1$ and $H_1\Delta_2$.  Second, on the open side
the local model should be a uniform expression, and not an arbitrary
product of thickened facets.  (This seems to be demanded on aesthetic and
logical grounds.  This author does not know if this will affect the final
outcome.)  The third and final point is that the relation between compact
and open is similar to that between $U^i$ and $U^j$, for $i+j=n$ and
$i\geq j$.  As noted, enlarging $U^i$ might reduce $U^j$ and hence reduce
$H^i$.  On the other hand, it might not change $U^j$, and can thus enlarge
$H^i$.  This phenomenom helps explain some of the subtle differences
between various similar-looking local-global groups.

We will now assume that whenever a system of compact local-global cycles is
defined, it will be paired with the space of open local-global cycles that
it determines (as with $U^i$ and $U^j$), to produce both the compact and
the open local-global groups.  (Recall that this is done by factoring out
the nullspaces.)

To complete this section we will exhibit, for $n\leq 6$, the appropriate
spaces of compact local-global cycles.  They all have the general form
$(\Delta_1 - \Delta_2)\eta$, where $\eta$ is a non-uniform expression.  As
already noted, for $n=3$ we take $\eta=\alpha\frown\beta\in T^2$.  This
gives, as promised $4+1=5=F_5$ homology groups.  (The $4$ comes from
$H_0$, $\ldots$, $H_4$.)  For $n=4$ we use $T^2$, $T^2\otimes S^1$ and
$T^3$ to obtain $5+3=8$ homology groups.

Now for $n=5$.  We use $T^2$, $T^2\otimes S^1$ and $T^2 \otimes S^2$. This
corresponds to a `string' of Betti numbers.  Because the Lefschetz element
$\omega$ is an element of $S^1$, it acts on this string of compact
local-global groups.  There will also be an adjoint map on the open
groups.  We also have $T^3$ and $T^3\otimes S^1$ as another string.  (This
group will be reconsidered shortly.)  And there is also $T^4$.  This gives
$6$ local-global groups organised into $3$ strings. Theere is also $H_0$,
$\ldots$, $H_5$.  As the general theory, as already noted, gives $13$
organised into $5$ strings, we have missed one of the local-global groups.

A certain amount of thought shows that the missing group consists of
compact local $2$-cycles.  Taking $\eta\in T^3$ is the way to produce such
a cycle, but this must be distinguished from the previous use of $T^3$. 
There, the $T^3$ represents a `one-dimensional family of local
$1$-cycles', and such will usually have a non-empty intersection against
the generic element of $S^1$, whereas such will always miss a local
$2$-cycle.  Thus, the `missing' group will be produced by $T^3 \cap
(S^1)^\perp$.  (By this is meant formal sums $\eta$ of $\alpha_1\frown
\alpha_2\frown \alpha_3$ such that $\eta \frown \xi$ is zero for any open
$\xi$ that is locally of the form $S^1 \otimes T^1$.  At this point it
make no sense to impose global conditions of $\xi$.)

The first use of $T^3$ is not, according to the local-global theory, quite
as it should be.  The cone on a simple $4$-polytope will in general have
non-trivial such $T^3$ items, but will not of course have any
`one-dimensional family of local $1$-cycles'.  The solution is to impose
an additional condition, besides consistency, on the open cycles that are
used.  THe condition is that if $\Delta_1$ and $\Delta_2$ differ only over
isolated points of $\Delta$, then the open local cycle $\xi$ should vanish
against $(\Delta_1-\Delta_2)\frown\eta$, for any $\eta$ in $T^3$.  A
similar condition should be imposed for $T^2$ when $n=4$, but in that case
it is vacuous.

For $n=6$ we will have all the $n=5$ items, both as themselves, and
multiplied by $S^1$.  In other words, there is $T^2$, $\ldots$,
$T^2\otimes S^3$ and $T^3$, $T^3\otimes S^1$, $T^3 \otimes S^2$ and
$T^4$, $T^4 \otimes S^1$.  The missing groups from $n=5$ produces
\[
    T^3 \cap (S^2)^\perp    \qquad
    (T^3 \otimes S^1) \cap (S^2)^\perp    
\]
for $n=6$.  (This last requires some thought.)  There will also be $T^5$
and $T^4 \cap (S^1)^\perp$.  Altogether this is $13$ local-global groups
organised into $6$ strings.  Conditions similar to those formulated in the
previous paragraph should again be imposed on open cycles.  Add to this the
$H_0$, $\ldots$, $H_6$ string to obtain $20$ groups in $7$ strings.  This
is short, by a string containing a single group.  This `missing' item will
be defined in the next section.

\section {Second order local-global homology}

In the previous section it was seen how suitable expressions of the form
$(\Delta_1 - \Delta_2)\frown\eta$ can be made to represent compact
local-global homology classes.  These, together with the uniform
expressions of \S3, were enough in dimension at most $5$ to produce the
$F_{n+2}$ homology groups that are demanded by the general combinatorial
structure of the local-global theory.  We also saw that for $n=6$, the
count was one short.  This section will define the missing group.  (The
previous sections dealt with order zero and order one respectively.)

First let $\Delta$ be the cone on a square.  This polytope has dimension
three, and in some sense is the first non-simple polytope.  On it there is
a non-trivial compact local-global cycle.  Now consider the product
$\Delta\times\Delta$.  This has dimension six.  On it one can form the
product $\eta$ of the local-global cycles on its factors.  Once this has
been suitably understood, it will provide an example of a second-order
local-global cycle.

It is natural to stratify $\Delta$ into the apex (the only non-simple
point) and the rest.  In the same way, $\Delta\times\Delta$ can be
stratified into $\{0\}\times\{0\}$, $\{0\}\times\Delta$,
$\Delta\times\{0\}$ and the rest.  Here `$0$' denotes the apex of
$\Delta$.  The cycle $\eta$ will in some sense be local to the `apex'
$\{0\}\times\{0\}$ of $\Delta\times\Delta$.  It represents a local-global
cycle that is local to $\{0\}\times\Delta$ (or to $\Delta\times\{0\}$)
that can further be made local to $\{0\}\times\{0\}$.  In fact, there will
be two such cycles, one for each factor in $\Delta\times\Delta$.  Although
equivalent when regarded as first-order local-global cycles, they will be
inequivalent when regarded as second-order such.

For simplicity, let $\Delta_a$ and $\Delta_b$ denote the two simple
polytopes formed from $\Delta$ by the process of $\epsilon$-variation. 
Now let $\Delta_{aa}$, $\ldots$, $\Delta_{bb}$ be the corresponding
resolutions of $\Delta\times\Delta$.  The alternating sum
\[
    \Delta_{{.}{.}} = \Delta_{aa} -  \Delta_{ab} 
                                - \Delta_{ba} +  \Delta_{bb}
\]
has the interesting property that above the whole of $\Delta\times\Delta$,
except for the apex $\{0\}\times\{0\}$, it is so to speak zero.  Consider,
for example, the cycle $\Delta_{{.}{.}}\frown\eta$, for any suitable
(non-uniform) $\eta$.  This will be a compact cycle that is, in some
sense, concentrated at the apex $\{0\}\times\{0\}$.  (Throughout this
section we will take $\eta$ to be as follows.  Each thickened facet
$\alpha$ on $\Delta$ determines thickened facets $\alpha\times\Delta$ and
$\Delta\times\alpha$ on $\Delta\times\Delta$.  We will have $\eta$ be
$(\psi\times\Delta)\frown(\Delta\times\psi)$ where $\psi$ is say $N\frown
E$ on $\Delta$, or any other non-uniform element of $T^2\Delta$.)

As presented in the previous paragraph, the expression
$\Delta_{{.}{.}}\frown\eta$ represents a first-order local-global cycle. 
Let us now present it as a second order such.  To do this we introduce
\[
    \Delta_{{.}a} = \Delta_{aa} - \Delta_{ba} \> , \qquad
    \Delta_{{.}b} = \Delta_{ab} - \Delta_{bb}
\]
and then consider both $\Delta_{{.}a}\frown\eta$ and
$\Delta_{{.}b}\frown\eta$.  Each of these expressions represents a
first-order local-global cycle on $\Delta\times\Delta$.  The difference
\[
    (\Delta_{{.}a} - \Delta_{{.}b}) \frown \eta
\]
is again a local-global cycle, but now it is local to the apex
$\{0\}\times\{0\}$.

The differences between the various local-global cycles that can be
constructed from $\Delta_{{.}{.}}$ and $\eta$ are subtle, and somewhat as
in the previous section they only become apparent once the dual open
theory has been defined.  For order zero `open' cycles, the basic model is
a (thickened) facet lying on $\Delta$.  For the first order theory, the
basic model is a facet passing through a vertex.  For the second order
theory, the basic model will be a facet that contains a flag.  Here, a
flag is for example a vertex lying on a $3$-face.  The `Betti numbers' are
in some clever way counting how many flags there are of each type.

Now once again consider the difference $\Delta_{{.}a} - \Delta_{{.}b}$. 
We do this not as the formal sum $\Delta_{{.}{.}}$, but as an expression
in its own right.  In other words, it can more exactly be thought of as
an ordered pair $(\Delta_{{.}a},\Delta_{{.}b})$, that is to be treated in
a particular way.  The quantity $\Delta_{{.}a}$ has a face of
$\Delta\times\Delta$ naturally associated to it, namely the face
$\{0\}\times\Delta$ along which the two components $\Delta_{aa}$ and
$\Delta_{ba}$ differ.  The same applies to $\Delta_{{.}b}$.  Now note that
$\Delta_{{.}a}$ and $\Delta_{{.}b}$ are so to speak concentrated over or
around the face $\{0\}\times\Delta$, and so it makes sense that they
should be paired with products of thickened facets that pass through this
face.

As already noted, this pairing will vanish over all of
$\Delta\times\Delta$, except for the apex $\{0\}\times\{0\}$.  One way to
record this fact is to observe that
$\Delta_{{.}{.}}\frown\eta\frown\alpha\frown\beta$
will always be zero, if it is known that $\alpha$ lies in $S^1$.  This is
true both as a local and a global statement.

For the open cycles, there is a detail missing.  Although the face
$\{0\}\times\Delta$ has entered the discussion, the apex
$\{0\}\times\{0\}$ has not.  It does so in the following way.  Let
$\alpha$ and $\beta$ be the thickened facets $\Delta\times N$ and
$N\times\Delta$ on $\Delta\times\Delta$.  (Any of the four triangular
facets of $\Delta$ could just as well have been chosen.)  Now consider the
object that is the flag 
$(\,\{0\}\times\{0\} \subset \{0\}\,)\times\Delta$,
together with the (open) expression $\alpha\frown\beta$.  Note that
expression has certain vanishing or incidence properties with respect to
the flag part of the object.  Both $\alpha$ and $\beta$ pass through
$\{0\}\times\{0\}$, while $\beta$ also contains $\{0\}\times\Delta$.  As
in \cite{bib.JF.CPLA}, these properties can be presented in an abstract
form.

This study of $\Delta\times\Delta$ has led to a compact cycle
$(\Delta_{{.}a}-\Delta_{{.}b}) \frown\eta$, and an open cycle (described in
the previous paragraph), that have a non-zero pairing with each other.  By
writing down the abstract properties that these cycles satisfy, a
definition of the second-order local-global cycles for $n=6$ will be
obtained.

At this point, it is well to stop.  The reader who is already familiar
with \cite{bib.JF.LGIH} and \cite{bib.JF.CPLA}, especially the former,
will appreciate that there are subtleties in the dimension of flags to be
used, that have yet to manifest themselves.  For the other readers, there
is no easy way to explain further, other than to appeal to the principles
already announced in \cite{bib.JF.LGIH}.

It has not been the goal of this section, to produce a complete and
rigorous definition of the higher-order local-global homology.  Rather, it
has been to develop the concepts in a fairly natural manner, up to the
point where all the major features have been exposed.  The goal has been
more to show the existence of a such an approach, than to exhibit it in a
formal and rigorous manner.

\section {Summary and conclusions}

This section discusses the following.  First, application to the
combinatorics of general convex polytopes.  Second, application to more
general algebraic varieties.  Third, application of the volume polynomial
approach to other situations, such as the Voronoi polytope of a positive
definite quadratic form.  Finally, some remarks are made on how the
conjectures implicit in this paper might be proved.

If $\Delta$ is a simple poytope, the known facts regarding
$H_\bullet\Delta$ imply numerical conditions on the face vector that are,
in addition to being necessary, are also sufficient for the existence of a
simple polytope with given flag vector.  These facts are generation by the
(thickened) facets, the ring structure, the strong Lefschetz theorem, and
a formula for the Betti numbers in terms of the face vector (and vice
versa).  The numerical conditions are implicit in the properties of
$H_\bullet\Delta$.  The proof of necessity requires both the strong
Lefschetz theorem and some result in monomial rings and the like.  The
proof of sufficiency is a matter of finding a suitable ingenious
construction, and showing that it gives a polytope with the required face
vector.

For general polytopes, results about the homology object $H_\bullet\Delta$
will again produce necessary conditions on this time the flag vector of
$\Delta$.  It is still true that $H_\bullet\Delta$ is generated, in some
sense, by the possibly thickened facets of $\Delta$.  This is of course
true for the expressions $\eta$, $\psi$ and so forth, whether uniform or
not.  It also seems to be true for the simple $\epsilon$-variations
$\Delta_i$. Each ordering of the facets of $\Delta$ will determine the
combinatorial type of such a $\Delta_i$.  Simply move the first facet
outward until it is in general position, then the second, then the third,
and so on.  Even if not all $\Delta_i$ arise in this way, perhaps a
`spanning set' will so arise.  The nature of any ring-like structure that
might exist (in the simple case this gives the `pseudo-power
inequalities') is not so obvious.  Proof of the strong Lefschetz theorem
and formulae for the Betti numbers are likely, in general, to be deep
results.  (It follows from Bayer's example \cite[\S?]{bib.JF.LGIH} that
such is unlikely to hold for all the local-global homology groups.  There
seems to be a strong analogy or connection between the
$\epsilon$-variation process and the construction of the `secondary
polytope' \cite{bib.LB-IG-BS.DMSP,bib.LB-BS.FP}.)

Nonetheless, it is true that results of this nature regarding
$H_\bullet\Delta$ will imply, as in the simple case, necessary numerical
conditions on the $f$-vector of $\Delta$.  Even if such results are only
conjectural, the result will be numerical conjectures on $f\Delta$.  This
would be an advance, for we are at present without even any plausible
conjectures for the conditions on $f\Delta$ in dimension greater than $3$.

Perhaps the best way to explore this problem is to focus on $n=4$.  This
problem is hard enough to be instructive, and easy enough to be
accessible.  Provided the conjectural structure of $H_\bullet\Delta$ can
be well understood in this case, it will be possible to extend the
existing construction in the simple case \cite{bib.LB-CL.SMC}, so that it
deals with this new situation.  Such is probably an appropriate starting
point for the study of the subtleties and complexities of local-global
intersection homology in higher dimensions.  The complications should
precisely satisfy the requirements of the proof process.

Now let $X$ be an irreducible projective algebraic variety.  First, let
$X$ be a Schubert variety, or something similar.  Such varieties have been
extensively studied.  The uniform expression approach taken to the
intersection homology of $\PDelta$ should apply with significant change to
this new situation.  In particular, the $\epsilon$-variation and study of
volume method should still be valid.  It may be that the concepts
introduced in this paper are related to results and methods already known
and used in this more complicated context.

Now let $X$ be any irreducible projective algebraic variety.  The homology
will no longer be generated by the `facets', and so the volume approach
will no longer give as much information as homology does.  An elliptic
curve is the simplest example of this.  The volume approach is however
very atractive, and it would be nice if for such a general $X$ there were
an $\epsilon$-variation (of something) that would together with the volume
analogue record all the information that homology does.  This can be
thought of as a problem for nonsingular varieties only.  One would also
wish for this theory to provide a `lifting' for each resolution $X_i\to
X$, similar to the decomposition theorem lifting that exists for $\Delta_i$
and $\Delta$.  Even without this, the decomposition theorem can be used to
define the concept of a uniform expression, and so in the same way produce
local-global homology groups.

The general method of $\epsilon$-variation and volume can be applied in
other situations, although with what success is not yet clear.  Here is an
example.  Let $Q$ be a positive definite quadratic form on say $\bfR^n$. 
Inside $\bfR^n$ is the integer lattice $\bfZ^n$.  The \emph{Voronoi
polytope} $\Delta_Q$ of $Q$ consists of all points of $\bfR^n$ that are at
least as close to the origin as they are to any other lattice point, where
the distance is measured using $Q$.  Now vary $Q$, say by adding a
quadratic form $\epsilon$ that is close to zero.  This will vary
$\Delta_Q$ to $\Delta_Q(\epsilon)$.  In this situation it is proper to
study not only the volume of $\Delta_Q(\epsilon)$, but also the area (or
length or whatever) of its faces of the various dimensions.  Provided the
same programme can be carried out, as has been done for simple polytopes,
subtle combinatorial inequalities on $\Delta_Q$, and thus $Q$, are likely
to result.

One might also wish to apply similar methods to, say, centrally symmetric
polytopes, or to arrangements of hyperplanes in an affine space, or even
to graphs and hypergraphs.  The beauty of the $\epsilon$-variation and
volume approach, when or if it works, is that it produces the definition
of homology out of straightforward geometric concepts, and the
requirements of consistency. It is as if it provides a nucleus or seed,
out of which an ingenious and appropriate homology theory might emerge. 
But for such to be useful, the Betti numbers should be linear functions of
the `flag vector'.

We close this paper by saying a few words about proofs.  The local-global
theory is of combinatorial interest only when it produces homology groups
whose dimension is a linear function of the flag vector.  The prototype
for the method of proof is in some sense McMullen's elementary and largely
geometric proof \cite{bib.McM.SP} of strong Lefschetz for simple convex
polytopes.  What happens here is that there is a package of properties,
that is preserved as both the geometric realization of a combinatorial
type (continuous change) and the combinatorial type itself (discrete
change) is varied, as a result of passing from one simple convex polytope
to another.  This package includes a strengthened version of the strong
Lefschetz theorem, namely the `Riemann-Hodge-Minkowski' inequalities. 
There is an induction built into the process.

Now consider the problem of proving that the order zero (or the usual
middle perversity intersection homology) part of the theory has the
predicted Betti numbers.  This is true for polytopes with rational
vertices, as a consequence of Deligne's proof of the Weil conjectures.  To
prove such a result for general polytopes, a package and induction similar
to that used by McMullen seems to be required.  The author believes that
the local-global intersection homology groups will provide the at present
unknown part of this package.  Again, this is a problem that can usefully
be investigated in the case $n=4$.

\end{document}